\documentclass[a4paper,12pt]{amsart}

\usepackage{amsfonts}
\usepackage{amssymb}
\usepackage{ifthen}
\usepackage{graphicx}
\usepackage{color}
\newtheorem{thm}{Theorem}
\newtheorem{cor}[equation]{Corollary}
\newtheorem{lem}{Lemma}
\newtheorem{prop}[equation]{Proposition}

\newtheorem{case}{Case}
\newtheorem{claim}{Claim}
\newtheorem{conj}[equation]{Conjecture}
\newtheorem{rem}{Remark}
\theoremstyle{definition}
\newtheorem{defn}{Definition}[section]
\newtheorem{example}{Example}[section]

\newtheorem{prob}[equation]{Problem}
\newtheorem{ques}[equation]{Question}

\newcounter {own}
\def\theown {\thesection       .\arabic{own}}

\newenvironment{pf}[1][]{%
 \vskip 3mm
 \noindent
 \ifthenelse{\equal{#1}{}}%
  {{\slshape Proof. }}%
  {{\slshape #1.} }%
 }%
{\qed\bigskip}

\newcounter{alphabet}
\newcounter{tmp}
\newenvironment{Thm}[1][]{\refstepcounter{alphabet}%
\bigskip%
\noindent%
{\bf Theorem \Alph{alphabet}}%
\ifthenelse{\equal{#1}{}}{}{ (#1)}%
{\bf .} \itshape}{\vskip 8pt}

\makeatletter
\newcommand{\Ref}[1]{\@ifundefined{r@#1}{}{\setcounter{tmp}{\ref{#1}}\Alph{tmp}}}
\makeatother

\newenvironment{Lem}[1][]{\refstepcounter{alphabet}%
\bigskip%
\noindent%
{\bf Lemma \Alph{alphabet}}%
{\bf .} \itshape}{\vskip 8pt}

\newcommand{\IR}{{\mathbb R}}

\newcommand{\ID}{{\mathbb D}}
\newcommand{\IB}{{\mathbb B}}

\def\be{\begin{equation}}
\def\ee{\end{equation}}

\newcommand{\bee}{\begin{enumerate}}
\newcommand{\eee}{\end{enumerate}}

\newcommand{\blem}{\begin{lem}}
\newcommand{\elem}{\end{lem}}
\newcommand{\bthm}{\begin{thm}}
\newcommand{\ethm}{\end{thm}}
\newcommand{\bcor}{\begin{cor}}
\newcommand{\ecor}{\end{cor}}
\newcommand{\beg}{\begin{examp}}
\newcommand{\eeg}{\end{examp}}
\newcommand{\begs}{\begin{examples}}
\newcommand{\eegs}{\end{examples}}
\newcommand{\bdefe}{\begin{defn}}
\newcommand{\edefe}{\end{defn}}
\newcommand{\bprob}{\begin{prob}}
\newcommand{\eprob}{\end{prob}}
\newcommand{\bques}{\begin{ques}}
\newcommand{\eques}{\end{ques}}
\newcommand{\bei}{\begin{itemize}}
\newcommand{\eei}{\end{itemize}}
\newcommand{\bcl}{\begin{claim}}
\newcommand{\ecl}{\end{claim}}
\newcommand{\bca}{\begin{case}}
\newcommand{\eca}{\end{case}}

\newcommand{\bcon}{\begin{conj}}
\newcommand{\econ}{\end{conj}}
\newcommand{\bcons}{\begin{conjs}}
\newcommand{\econs}{\end{conjs}}
\newcommand{\bprop}{\begin{prop}}
\newcommand{\eprop}{\end{prop}}
\newcommand{\br}{\begin{rem}}
\newcommand{\er}{\end{rem}}
\newcommand{\brs}{\begin{rems}}
\newcommand{\ers}{\end{rems}}
\newcommand{\bo}{\begin{obser}}
\newcommand{\eo}{\end{obser}}
\newcommand{\bos}{\begin{obsers}}
\newcommand{\eos}{\end{obsers}}
\newcommand{\bpf}{\begin{pf}}
\newcommand{\epf}{\end{pf}}
\newcommand{\ba}{\begin{array}}
\newcommand{\ea}{\end{array}}
\newcommand{\beq}{\begin{eqnarray}}
\newcommand{\beqq}{\begin{eqnarray*}}
\newcommand{\eeq}{\end{eqnarray}}
\newcommand{\eeqq}{\end{eqnarray*}}

\newcommand{\ds}{\displaystyle}

\newcounter{minutes}
\divide\time by 60
\newcounter{hours}
\multiply\time by 60 \addtocounter{minutes}{-\time}

\begin{document}
\bibliographystyle{amsplain}
\title [] {Lipschitz type spaces and Landau-Bloch type theorems for harmonic functions and Poisson equations }

\def\thefootnote{}
\footnotetext{ \texttt{\tiny File:~\jobname .tex,
          printed: \number\day-\number\month-\number\year,
          \thehours.\ifnum\theminutes<10{0}\fi\theminutes}
} \makeatletter\def\thefootnote{\@arabic\c@footnote}\makeatother

\author{Sh. Chen}
\address{Sh. Chen,
Department of Mathematics and Computational Science, Hengyang Normal
University, Hengyang, Hunan 421008, People's Republic of China.}
\email{mathechen@126.com}

\author{M. Mateljevi\'c}
\address{M. Mateljevi\'c,
University of Belgrade, Faculty of Mathematics
Studentski Trg 16, 11000 Belgrade, Serbia.} \email{miodrag@matf.bg.ac.rs}

\author{S. Ponnusamy$^\dagger $
}
\address{S. Ponnusamy,
Indian Statistical Institute (ISI), Chennai Centre, SETS (Society
for Electronic Transactions and security), MGR Knowledge City, CIT
Campus, Taramani, Chennai 600 113, India. }
\email{samy@isichennai.res.in, samy@iitm.ac.in}

\author{X. Wang
}
\address{X. Wang, Department of Mathematics,
Hunan Normal University, Changsha, Hunan 410081, People's Republic
of China.} \email{mathshida@gmail.com}

\subjclass[2000]{Primary: 47B33, 31C25}
\keywords{Schwarz's lemma,  harmonic function, Lipschitz type space,
Poisson equation.
\\
$
^\dagger$ {\tt This author is on leave from the Department of Mathematics,
Indian Institute of Technology Madras, Chennai-600 036, India.}\\
}

\begin{abstract}
In this paper, we investigate some properties on harmonic functions
and solutions to Poisson equations. First, we will discuss the
Lipschitz type spaces on harmonic functions. Secondly, we establish
the Schwarz-Pick lemma for harmonic functions in the unit ball
$\IB^n$ of $\IR^n$, and then we apply it to obtain a Bloch theorem
for harmonic functions in Hardy spaces. At last, we use a normal
family argument to extend the Landau-Bloch type theorem  to
functions which are solutions to Poisson equations.
\end{abstract}



\maketitle \pagestyle{myheadings} \markboth{SH. Chen,  M.
Mateljevi\'c, S. Ponnusamy and  X. Wang} {Lipschitz type spaces and
Landau-Bloch type theorems }

\section{Introduction and main results}\label{csw-sec1}
Let $\mathbb{R}^{n}
$ denote the usual real vector space of dimension $n$, where
$n\geq2$ is a positive integer. Sometimes it is convenient to
identify each point $x=(x_{1},\ldots,x_{n})\in\mathbb{R}^{n}$ with
an $n\times 1$ column matrix so that
$$x=\left(\begin{array}{cccc}
x_{1}   \\
\vdots \\
 x_{n}
\end{array}\right).
$$
For $a=(a_{1},\ldots,a_{n})$ and $x\in\mathbb{R}^{n}$, we define the Euclidean inner product
$\langle \cdot ,\cdot \rangle$ by
$$\langle x,a\rangle=x_{1}a_{1}+\cdots+x_{n}a_{n}
$$
so that the Euclidean length of $x$ is defined by
$$|x|=\langle x,x\rangle^{1/2}=(|x_{1}|^{2}+\cdots+|x_{n}|^{2})^{1/2}.
$$
Denote a ball in $\mathbb{R}^{n}$
with center $x'$ and radius $r$ by
$$\mathbb{B}^{n}(x',r)=\{x\in\mathbb{R}^{n}:\, |x-x'|<r\}.
$$
In particular, $\mathbb{B}^{n}$ denotes the unit ball
$\mathbb{B}^{n}(0,1)$. Set $\mathbb{D}=\mathbb{B}^2$, the open unit
disk in the complex plane $\mathbb{C}$.

A function $f$ of an open subset $\Omega\subset\mathbb{R}^{n}$ into
$\mathbb{R}$ is called a {\it harmonic function} if $\Delta f=0$,
where $\Delta$ represents the $n$-dimensional Laplacian operator
$$\Delta=\sum_{k=1}^{n}\frac{\partial^{2}}{\partial x_{k}^{2}}.
$$

In this paper, we use $C$ to denote the various positive constants,
whose value may change from one occurrence to the next.

A continuous increasing function $\omega:\, [0,\infty)\rightarrow
[0,\infty)$ with $\omega(0)=0$ is called a {\it majorant} if
$\omega(t)/t$ is non-increasing for $t>0$. Given a subset $\Omega$
of $\mathbb{R}^{n}$, a function $f:\, \Omega\rightarrow
\mathbb{R}^{m}~(m\geq1)$ is said to belong to the {\it Lipschitz
space $\Lambda_{\omega}(\Omega)$} if there is a positive constant
$C$ such that \be\label{eq1x} |f(x)-f(y)|\leq C\omega(|x-y|) ~\mbox{
for all $x,\ y\in\Omega.$} \ee For $\delta_{0}>0$, let
\be\label{eq2x} \int_{0}^{\delta}\frac{\omega(t)}{t}\,dt\leq
C\cdot\omega(\delta),\ 0<\delta<\delta_{0} \ee and \be\label{eq3x}
\delta\int_{\delta}^{\infty}\frac{\omega(t)}{t^{2}}\,dt\leq
C\cdot\omega(\delta),\ 0<\delta<\delta_{0}, \ee where $\omega$ is a
majorant. A majorant $\omega$ is said to be {\it regular} if it
satisfies the conditions (\ref{eq2x}) and (\ref{eq3x}) (see
\cite{CRW,D,D1,P,Pav1,Pav2}).

Let $\Omega$ be a domain in $\mathbb{R}^{n}$ with non-empty
boundary. We use $d_{\Omega}(x)$ to denote the Euclidean distance
from $x$ to the boundary $\partial \Omega$ of $\Omega$. In
particular, we always use $d(x)$ to denote the Euclidean distance
from $x$ to the boundary of $\mathbb{B}^{n}.$

A proper subdomain $G$ of $\mathbb{R}^{n}$  is said to be {\it
$\Lambda_{\omega}$-extension} if
$\Lambda_{\omega}(G)=\mbox{loc}\Lambda_{\omega}(G)$, where
$\mbox{loc}\Lambda_{\omega}(G)$ denotes the set of all functions
$f:\, G\rightarrow \mathbb{R}^{m}$ satisfying  (\ref{eq1x}) with a
fixed positive constant $C$, whenever $x\in G$ and $y\in G$ such
that $|x-y|<\frac{1}{2}d_{G}(x)$.  Obviously,  $\mathbb{B}^{n}$ is a
$\Lambda_{\omega}$-extension domain.

In \cite{L}, the author proved that $G$ is a
$\Lambda_{\omega}$-extension domain if and only if each pair of
points $x,y\in G$ can be joined by a rectifiable curve
$\gamma\subset G$ satisfying \be\label{eq1.0}
\int_{\gamma}\frac{\omega(d_{G}(\zeta))}{d_{G}(\zeta)}\,ds(\zeta)
\leq C\omega(|x-y|) \ee with some fixed positive constant
$C=C(G,\omega)$, where $ds$ stands for the arc length measure on
$\gamma$.  Furthermore, Lappalainen \cite[Theorem 4.12]{L} proved
that $\Lambda_{\omega}$-extension domains  exist only for majorants
$\omega$ satisfying  (\ref{eq2x}). See \cite{D1,GM,KW,L} for more
details on $\Lambda_{\omega}$-extension domains.

Krantz \cite{Kr} proved a Hardy-Littlewood type theorem for harmonic
functions in the unit ball with respect to the majorant
$\omega(t)=\omega_{\alpha}(t)=t^{\alpha}~(0<\alpha\leq1)$ as
follows.

\begin{Thm}{\rm (\cite[Theorem 15.8]{Kr})}\label{ThmA}
Let $f$ be a harmonic function from $\mathbb{B}^{n}$ into
$\mathbb{R}$ and $0<\alpha\leq1$. Then $f$ satisfies $$|\nabla
f(x)|\leq C\frac{\omega_{\alpha}\big(d(x)\big)}{d(x)}~\mbox{for
any}~x\in\mathbb{B}^{n}$$ if and only if
$$|f(x)-f(y)|\leq C\omega_{\alpha}(|x-y|)~\mbox{for
any}~x,y\in\mathbb{B}^{n},$$  where $\nabla f$ denotes the gradient
of $f$.
\end{Thm}

For the extensive discussions on this topic, see
\cite{Ai,ABM,MVM,ABN,CPW5}. We generalize Theorem \Ref{ThmA} to the
following form.

\begin{thm}\label{thm-1}
Let $\omega$ be a majorant satisfying  {\rm (\ref{eq2x})}, $\Omega$
be a $\Lambda_{\omega}$-extension domain in $\mathbb{R}^{n}$ and $f$
be a harmonic function from $\Omega$ into $\mathbb{R}$. Then $f\in
\Lambda_{\omega}(\Omega)$  if and only if
$$|\nabla f(x)|\leq
C\frac{\omega\big(d_{\Omega}(x)\big)}{d_{\Omega}(x)}~\mbox{for
any}~x\in\Omega.$$
\end{thm}

In \cite{Ho}, Holland-Walsh  obtained the following result. For the
extensive studies on this topic, see \cite{CPW2,Pav,Re}.

\begin{Thm}{\rm (\cite[Theorem 3]{Ho})}\label{ThmA2}
Let  $\mathcal{B}$   denote all analytic functions  in $\mathbb{D}$
which form a complex Banach space with the norm
$$\|f\|_{\mathcal{B}}=|f(0)|+\sup_{z\in\mathbb{D}}\{(1-|z|^{2})|f'(z)|\}<\infty.
$$
Then $f\in\mathcal{B}$ if and only if
$$\sup_{z,w\in\mathbb{D},z\neq w}\left\{
\frac{\sqrt{(1-|z|^{2})(1-|w|^{2})}|f(z)-f(w)|}{|z-w|}\right\}<\infty.
$$
\end{Thm}

In \cite{Pav}, Pavlovi\'c generalized Theorem \Ref{ThmA2} into the
following form.

\begin{Thm}{\rm (\cite[Theorem 2]{Pav})}\label{ThmA3}
Let $ \mathcal{C}^{1}(\mathbb{B}^{n})$ be the class of all one order
continuous differentiable functions from $\mathbb{B}^{n}$ into
$\mathbb{R}.$ Let $\mathcal{B}_{\mathcal{C}^{1}}$ denote all
 $f\in \mathcal{C}^{1}(\mathbb{B}^{n})$ which form a
Banach space with the norm
$$\|f\|_{\mathcal{B}_{\mathcal{C}^{1}}}=|f(0)|+\sup_{z\in\mathbb{D}}\{(1-|x|^{2})|\nabla
f(x)|\}<\infty.
$$
Then $f\in\mathcal{B}_{\mathcal{C}^{1}}$ if and only if
$$\sup_{x,y\in\mathbb{B}^{n},x\neq y}\left\{
\frac{\sqrt{(1-|x|^{2})(1-|y|^{2})}|f(x)-f(y)|}{|x-y|}\right\}<\infty.
$$
\end{Thm}


By using  a different proof methods,
 we will prove a more general result as follows which is
a generalization of Theorems \Ref{ThmA2} and \Ref{ThmA3}.

\begin{thm}\label{thm-CPW1}
Let $f\in\mathcal{C}^{1}(\mathbb{B}^{n})$ and $\omega$ be a
majorant. Then for any $x\in\mathbb{B}^{n}$,
$$|\nabla f(x)|\leq C\omega\left(\frac{1}{d(x)}\right)
$$ if and only if
for any $x, y\in\mathbb{R}^{n}$ with $x\neq y$,
$$\frac{|f(x)-f(y)|}{|x-y|}\leq
C\omega\left(\frac{1}{\sqrt{d(x)d(y)}}\right).
$$

\end{thm}

Dyakonov \cite{D} discussed the relationship between the Lipschitz
space and the bounded mean oscillation on holomorphic functions in
$\mathbb{D}$, and obtained the following result.

\begin{Thm}{\rm (\cite[Theorem 1]{D})}\label{ThmDy}
Suppose that $f$ is a holomorphic function in $\mathbb{D}$ which is
continuous up to the boundary of $\mathbb{D}$. If $\omega$ and
$\omega^{2}$ are regular majorants,  then
$$f\in \Lambda_{\omega}(\mathbb{D})\Longleftrightarrow \mbox{{\rm P}}_{|f|^{2}}(z)-|f(z)|^{2}\leq M\omega^{2}(d(z)),
$$
where
$$\mbox{{\rm P}}_{|f|^{2}}(z)=\frac{1}{2\pi}\int_{0}^{2\pi}\frac{1-|z|^{2}}{|z-e^{i\theta}|^{2}}|f(e^{i\theta})|^{2}\,d\theta.
$$
\end{Thm}

In particular, for harmonic functions, we get the following result
which is analogous to Theorems \Ref{ThmA2}  and \Ref{ThmDy}. For
some related topics on complex-valued functions, we refer to
\cite{CPVW,CRW}.

\begin{thm}\label{thm-CPW}
Let $f\in\mathcal{C}^{1}(\mathbb{B}^{n})$ be a harmonic  and
$\omega$ be a majorant. Then the following are equivalent:

\item{{{\rm(a)}}}~for any $x\in\mathbb{B}^{n}$, $$|\nabla f(x)|\leq C\omega\left(\frac{1}{d(x)}\right);
$$

\item{{{\rm(b)}}}~ for any $x, y\in\mathbb{R}^{n}$ with $x\neq y$,
$$\frac{|f(x)-f(y)|}{|x-y|}\leq
C\omega\left(\frac{1}{\sqrt{d(x)d(y)}}\right);
$$

\item{{{\rm(c)}}}~for any  $r\in(0,d(x)]$,
$$\frac{1}{|\mathbb{B}^{n}(x,r)|}\int_{\mathbb{B}^{n}(x,r)}|f(\zeta)-f(x)|dV(\zeta)\leq Cr\omega\Big(\frac{1}{r}\Big),$$
where $dV$ denotes the Lebesgue volume measure in $\mathbb{B}^{n}$.
\end{thm}

For a {\it vector-valued and real harmonic
function}  $f=(f_{1},\ldots,f_{n})$ from $\mathbb{B}^{n}$ into  $\mathbb{R}^n$ (i.e.  for
each $i\in\{1,2,\ldots,n\}$, $f_{i}:\,\, \mathbb{B}^{n}\rightarrow \mathbb{R}$ is harmonic),
we denote the Jacobian of $f$ by $J_{f}$, i.e.,
$$J_{f}=\det \left ( \frac{\partial f_i}{\partial x_{j}}\right )_{n\times n},
$$
where $j\in\{1,2,\ldots,n\}.$ Let $\mathcal{H}(\mathbb{B}^{n}, \mathbb{R}^{n})$ be the set of all
real harmonic functions $f$ from $\mathbb{B}^{n}$ into
$\mathbb{R}^{n}$. Also, for $p\in(0,\infty)$, let
$\mathcal{H}^{p}(\mathbb{B}^{n}, \mathbb{R}^{n})$ denote the
harmonic Hardy class consisting of all functions
$f\in\mathcal{H}(\mathbb{B}^{n}, \mathbb{R}^{n})$ such that
$$\|f\|_{p}=\sup_{0<r<1}M_{p}(f,r) <\infty, \quad
M_{p}^{p}(f,r):=\int_{\partial\mathbb{B}^{n}}|f(r\zeta)|^{p}\, d\sigma(\zeta),
$$
where $d\sigma$ is the normalized surface measure on
$\partial\mathbb{B}^{n}$ (see \cite{ABR}).





One of the long standing open problems in geometric function theory
is to determine the precise value of the univalent Landau-Bloch
constant for analytic functions of $\ID$. It has attracted much
attention, see \cite{LM, Mi1,Mi2,M-89} and references therein. For
general holomorphic mappings of more than one complex variable, no
univalent Landau-Bloch constant exists (cf. \cite{ W}). In order to
obtain some analogous results of univalent Landau-Bloch constant for
functions with several complex variables, it is necessary to
restrict the class of mappings considered, see
\cite{CPW3,CPW,FG,LX,M3,T,W}.

In \cite{HG1}, the authors discussed the Schwarz-Pick Lemma and the
 Landau-Bloch type theorems for bounded pluriharmonic mappings. It is known that pluriharmonic
 mappings are special vector-valued harmonic
functions.
   By using a different approach, as our last aim, we will
 establish the Schwarz-Pick Lemma and obtain a univalent Landau-Bloch constant for vector-valued harmonic
functions in the Hardy spaces.  Since all bounded vector-valued
harmonic functions belong to the harmonic Hardy classes,  we see
that our result (Theorem \ref{thm1}) is a generalization of
\cite[Theorem 5]{HG1}.

\begin{thm}\label{thm1}
Suppose that $f\in\mathcal{H}^{p}(\mathbb{B}^{n},\mathbb{R}^{n})$
satisfies  $J_{f} (0)-1=|f(0)|=0,$ where $p\geq1$ and $n\geq3$. Then
$f(\mathbb{B}^{n})$ contains a univalent ball $\mathbb{B}^{n}(0,R)$,
where
$$R\geq\max_{0<r<1}\varphi(r),
$$
where
$$\varphi(r)=\frac{1}{2[nK(r)]^{2n-2}M(r)[(1+\sqrt{2})^{n-1}+\sqrt{2}-1]},
$$
$$K(r) = 2^{1/p}\|f\|_{p}/[r(1-r)^{(n-1)/p}]~\mbox{and}~M(r)=K(r)[(3+\sqrt{3})n+2\sqrt{2}].
$$

\end{thm}

We remark that, as $\lim_{r\rightarrow0+}\varphi(r)=\lim_{r\rightarrow1-}\varphi(r)=0,
$
the maximum of $\varphi(r)$ in Theorem \ref{thm1} does exist.

The following result easily follows from Theorem \ref{thm1}.

\begin{thm}\label{thm5-x}
Let $f\in\mathcal{H}(\mathbb{B}^{n}, \mathbb{R}^{n})$ with  $J_{f}
(0)-1=|f(0)|=0$ and $|f(x)|<M$ for $x\in\mathbb{B}^{n}$. Then $f$ is
univalent in $\mathbb{B}^{n}(0,\rho_{0})$ and
$f(\mathbb{B}^{n}(0,\rho_{0}))$ contains a univalent ball
$\mathbb{B}^{n}(0,R_{0})$, where
$$\rho_{0}=\frac{1}{n^{n-1}M^{n}[(3+\sqrt{2})n+2\sqrt{2}][(1+\sqrt{2})^{n-1}+\sqrt{2}-1]}
~\mbox{ and }~R_{0}=\frac{\rho_{0}}{2(nM)^{n-1}}.
$$
\end{thm}

We will extend Theorem \ref{thm5-x} to a general case. Let us give
some preparations before we present our next result.

Let $f:~\overline{\Omega}\to \mathbb{R}^n$ be a differentiable
mapping and $p$ be a regular value of $f$, where $p\notin
f(\partial\Omega)$ and $\Omega\subset \mathbb{R}^n$ is a bounded
domain. Then the degree $\deg(f,\Omega,p)$ is defined by the formula
$$\deg(f,\Omega,p):=\sum_{y\in f^{-1}(p)}\mbox{ sign} \big(\det J_{f} (y)\big).
$$
The $\deg(f,\Omega,p)$ satisfies the following properties (cf. \cite{RR,V}):\\
\begin{enumerate}
\item[(I)]    If $\deg(f,\overline{\Omega},p)\neq 0$, then there
exists an $x\in\Omega$  such that $f(x)=p$.
\item[(II)]\label{(II)} If $D$ is a domain with $\overline{D}\subset\Omega$ and $p\in
\mathbb{R}^ n \backslash f(\partial D)$,    then the degree
$\deg(f,D,p)$ is a constant.
\end{enumerate}

Let $D\subset\mathbb{R}^{n}$  be a domain and $f$ be a real function
from $D$ into $\mathbb{R}^{n}$. If the H\"{o}lder coefficient
$$\|f\|_{\alpha,D}= \sup_{x,y\in D, x\neq y}  \frac{|f(x)- f(y)|}{|x-y|^\alpha}
$$
is finite, then the function $f$ is said to be (uniformly)
H\"{o}lder continuous with exponent $\alpha$ in $D$, where $0<
\alpha\leq 1.$ In this case, the H\"{o}lder coefficient serves as a
seminorm. If the H\"older coefficient is merely bounded on compact
subsets of $D$, then the function $f$ is said to be locally H\"older
continuous with exponent  $\alpha$ in $D$. We denote by
$C^\alpha(D,\mathbb{R}^{n})$ the space consist of all   locally
H\"{o}lder continuous functions $f$ from  $D$ into $\mathbb{R}^{n}$
with exponent $\alpha$ (cf. \cite{gt,Kr}).


Let $\mathcal{PE}_{f}$ denote the class of functions $u$ satisfying
the Poisson equation $\Delta u=f $ with $J_{u} (0)-1=|u(0)|=0$,
where  $u\in C^2(\mathbb{B}^{n})$, i.e.,   twice continuously
differentiable function in $\mathbb{B}^{n}$, and $f\in
C^\alpha(\mathbb{B}^{n}, \mathbb{R}^{n})$ with the constants
$\alpha\in(0,1)$ and $\|f\|_{\alpha,\mathbb{B}^{n}}<\infty$. We use
$\mathcal{PE}_{f}^M$ to denote  the family of all functions $u$
satisfying $u \in \mathcal{PE}_{f}$ with $|u(x)|\leq M$ for
$x\in\mathbb{B}^{n}$, where $M$ is a positive constant. Obviously,
all bounded harmonic functions belong to $\mathcal{PE}_{f}^M$.

\begin{thm}\label{thm6-x}
Let $u\in\mathcal{PE}_{f}^M$. Then  there is a positive constant
$c_0$ depending only on $M$, $\|f\|_{\alpha,\mathbb{B}^{n}}$ and $n$
such that $\mathbb{B}^{n}(0,c_0)\subset u(\mathbb{B}^{n})$.
\end{thm}

In fact, the bounded condition in Theorem \ref{thm6-x} is necessary.
The following example shows that there is no  Landau-Bloch Theorem
for functions $u\in\mathcal{PE}_{f}$ without the bounded condition.

\begin{example}
For $k\in\{1,2,\ldots\}$ and $x\in\mathbb{B}^{n}$, let  $u_k(x)=
(kx_1, x_2/k,x_3,\ldots,x_n)$. Then $u_k$ are harmonic  and $J_{u_k}
(0)-1=|u_k (0)|=0$.

 This  example  tells us that if    $u: \mathbb{B}^{n}\rightarrow  \mathbb{R}^{n}$
is a  harmonic   function on the unit ball  with  $J_{u}
(0)-1=|u(0)|=0$,  then   there is no an absolute constant $s>0$ such
that   $\mathbb{B}^{n}(0,s) $ belongs to   $u(\mathbb{B}^{n})$. Thus
the Theorem \ref{thm6-x}  does not hold  for $u\in\mathcal{PE}_{f}$.
\end{example}


The proofs of Theorems \ref{thm-1}, \ref{thm-CPW1} and \ref{thm-CPW}
will be given in Section \ref{csw-sec2}. We will show Theorems
\ref{thm1} and \ref{thm6-x}
 in the last part of this paper.

\section{Lipschitz type spaces on harmonic functions}\label{csw-sec2}

\subsection*{Proof of Theorem \ref{thm-1}}
We first prove the sufficiency. Since $\Omega$ is a
$\Lambda_{\omega}$-extension domain in $\mathbb{R}^{n}$, we see that
for any $x,y\in\Omega$, by using (\ref{eq1.0}), there is a
rectifiable curve $\gamma\subset\Omega$ joining $x$ to $y$ such that

\begin{eqnarray*}
|f(x)-f(y)|&\leq&\int_{\gamma}|\nabla f(\zeta)|ds(\zeta)\\
&\leq&C\int_{\gamma}\frac{\omega\big(d_{\Omega}(\zeta)\big)}{d_{\Omega}(\zeta)}ds(\zeta)\\
&\leq&C\omega(|x-y|).
\end{eqnarray*}

Now we come to prove the necessity. Let
$x=(x_{1},\ldots,x_{n})\in\Omega$ and $r=d_{\Omega}(x)/2$. 
For all $y\in\mathbb{B}^{n}(x,r)$, using Poisson formula, we get
$$f(y)=\int_{\partial\mathbb{B}^{n}}\mbox{P}(y,\zeta)f(r\zeta+x)d\sigma(\zeta),
$$
where $\zeta=(\zeta_{1},\ldots,\zeta_{n})\in\partial\mathbb{B}^{n}$ and
$$\mbox{P}(y,\zeta)=\frac{r^{2}-|y-x|^{2}}{|y-x-r\zeta|^{2}} 
.$$
By elementary calculations, for each $k\in\{1,2,\ldots,n\}$, we have
$$\frac{\partial \mbox{P}(y,\zeta)}{\partial
y_{k}}=-2
\frac{\big[(y_{k}-x_{k})|y-x-r\zeta|^{2}+(r^{2}-|y-x|^{2})(y_{k}-r\zeta_{k}-x_{k})\big]}{|y-x-r\zeta|^{4}}.$$
Then  for all $y\in\mathbb{B}^{n}(x,r/2)$,

\begin{eqnarray*}
\left|\frac{\partial \mbox{P}(y,\zeta)}{\partial y_{k}}\right|
&\leq&2\frac{\big[|y_{k}-x_{k}||y-x-r\zeta|^{2}+(r^{2}-|y-x|^{2})|y_{k}-r\zeta_{k}-x_{k}|\big]}{|y-x-r\zeta|^{4}}\\
&\leq&2\frac{\left[\frac{r}{2}\Big(\frac{3r}{2}\Big)^{2}+r^{2}\Big(\frac{3r}{2}\Big)\right]}{\Big(r-\frac{r}{2}\Big)^{4}}
= \frac{84}{r},
\end{eqnarray*}

\noindent which implies that

\begin{eqnarray*}
|\nabla
f(y)|&=&\left[\sum_{k=1}^{n}f^{2}_{y_{k}}(y)\right]^{\frac{1}{2}}\\
&=&\Big\{\sum_{k=1}^{n}\Big(\Big|\int_{\partial\mathbb{B}^{n}}\frac{\partial
\mbox{P}(y,\zeta)}{\partial
y_{k}}(f(r\zeta+x)-f(x))d\sigma(\zeta)\Big|^{2}\Big\}^{\frac{1}{2}}\\
&\leq&\sum_{k=1}^{n}\Big|\int_{\partial\mathbb{B}^{n}}\frac{\partial\mbox{P}(y,\zeta)}{\partial
y_{k}}(f(r\zeta+x)-f(x))d\sigma(\zeta)\Big|\\
&\leq&\sum_{k=1}^{n}\int_{\partial\mathbb{B}^{n}}\left|\frac{\partial\mbox{P}(y,\zeta)}{\partial
y_{k}}\right|\big|f(r\zeta+x)-f(x)\big|d\sigma(\zeta)\\
&\leq&\sqrt{n}\int_{\partial\mathbb{B}^{n}}\left|\nabla
\mbox{P}(y,\zeta)\right|\big|f(r\zeta+x)-f(x)\big|d\sigma(\zeta)\\
&\leq&\frac{84n}{r}\int_{\partial\mathbb{B}^{n}}\big|f(r\zeta+x)-f(x)\big|d\sigma(\zeta)\\
&\leq&\frac{84nC\omega(r)}{r}\\
&\leq&168nC\frac{\omega\big(d_{\Omega}(x)\big)}{d_{\Omega}(x)}.
\end{eqnarray*}
If we take $y=x$, then we get the desired result. The proof of this
theorem is complete. \qed

\subsection*{Proof of Theorem \ref{thm-CPW1}}

We first prove the necessity. For any $x, y\in\mathbb{B}^{n}$ with
$x\neq y$, let $\varphi(t)=xt+(1-t)y$, where $t\in[0,1]$. Since
$|\varphi(t)|\leq t|x|+(1-t)|y|$, we see that
\be\label{eq-y1}1-|\varphi(t)|\geq1-t|x|-|y|+t|y| \geq1-t+|y|(t-1)
=(1-t)d(y)
\ee
and
\be\label{eq-y2}1-|\varphi(t)|\geq1-t|x|-|y|+t|y|
=1-t|x|-|y|(1-t)\geq1-t|x|-(1-t)=td(x).
\ee
By (\ref{eq-y1}) and (\ref{eq-y2}), we get
$$\left(1-|\varphi(t)|\right)^{2}\geq (1-t)t d(x)d(y),
$$
which implies
\be\label{eq-t20}
\frac{1}{1-|\varphi(t)|}\leq\frac{1}{\sqrt{(1-t)t d(x)d(y)}}.
\ee

For $t>0,$ by the monotonicity of $\omega(t)/t$, we know that
\be\label{eq-p1}\omega(\lambda t)\leq\lambda\omega(t),\ee where
$\lambda\geq1$.

By (\ref{eq-t20}) and (\ref{eq-p1}), for any $x, y\in\mathbb{B}^{n}$
with $x\neq y$, we have
\begin{eqnarray*}
|f(x)-f(y)|&=&\left|\int_{0}^{1}\frac{df}{dt}(\varphi(t))dt\right|\\
&\leq&\sqrt{n}|x-y|\int_{0}^{1}|\nabla f(\varphi(t))|dt\\
&\leq&\sqrt{n}|x-y|\int_{0}^{1}\frac{|\nabla
f(\varphi(t))|}{\omega\left(\frac{1}{1-|\varphi(t)|}\right)}
\omega\left(\frac{1}{1-|\varphi(t)|}\right)dt\\
&\leq&C\sqrt{n}|x-y|\int_{0}^{1}\omega\left(\frac{1}{1-|\varphi(t)|}\right)dt\\
&\leq&C\sqrt{n}|x-y|\int_{0}^{1}\omega\left(\frac{1}{\sqrt{(1-t)t
d(x)d(y)}}\right)dt\\
&\leq&C\sqrt{n}|x-y|\omega\left(\frac{1}{\sqrt{d(x)d(y)}}\right)\int_{0}^{1}\frac{1}{\sqrt{(1-t)t}}dt\\
&=&C\sqrt{n}|x-y|\omega\left(\frac{1}{\sqrt{d(x)d(y)}}\right)\int_{0}^{\frac{\pi}{2}}\frac{2\sin\theta\cos\theta}
{\sqrt{\sin^{2}\theta\cos^{2}\theta}}d\theta\\
&=&C\pi\sqrt{n}|x-y|\omega\left(\frac{1}{\sqrt{d(x)d(y)}}\right),
\end{eqnarray*}
which gives
$$\frac{|f(x)-f(y)|}{|x-y|}\leq \pi
C\sqrt{n}\omega\left(\frac{1}{\sqrt{d(x)d(y)}}\right).
$$

Now we prove the sufficiency part.  For any $x,
y\in\mathbb{B}^{n}$ with $x\neq y$, since
$$|\nabla f(x)|=\lim\sup_{y\rightarrow
x}\frac{|f(y)-f(x)|}{|y-x|},$$ we see that $$\lim\sup_{y\rightarrow
x}\frac{|f(y)-f(x)|}{|y-x|}=|\nabla f(x)|\leq
C\lim\sup_{y\rightarrow
x}\omega\left(\frac{1}{\sqrt{d(x)d(y)}}\right)=C\omega\left(\frac{1}{d(x)}\right).$$
The proof of this theorem is complete. \qed

\vspace{6pt}

Using  arguments  similar to those in the proof of \cite[Lemma
2.5]{MV}, we have the following lemma and so, we omit its proof.

\begin{lem}\label{lem-g1}
Suppose that $f:\
\overline{\mathbb{B}}^{n}(a,r)\rightarrow\mathbb{R}$ is a continuous
function in $\overline{\mathbb{B}}^{n}(a,r)$ and harmonic  in
$\mathbb{B}^{n}(a,r)$. Then
$$|\nabla f(a)|\leq
\frac{n\sqrt{n}}{r}\int_{\partial\mathbb{B}^{n}}|f(a+r\zeta)-f(a)|d\sigma(\zeta).$$
\end{lem}

\subsection*{Proof of Theorem \ref{thm-CPW}}
(a)$\Longleftrightarrow$(b) easily follows from Theorem
\ref{thm-CPW1}. We only need to prove (a)$\Longleftrightarrow$(c).
We first prove  (a)$\Longrightarrow$(c). For any
$x=(x_{1},\ldots,x_{n}),\ y=(y_{1},\ldots,y_{n})\in\mathbb{B}^{n}$
and $t\in[0,1]$, we have
$$d\big(x+t(y-x)\big)\geq d(x)-t|y-x|.$$ Suppose that
$d(x)-t|y-x|>0$. Then

\begin{eqnarray*}
|f(y)-f(x)|&=&\Big|\int_{0}^{1}\frac{df}{dt}(\varsigma)dt\Big|\\
&=&\left|\sum_{k=1}^{n}(y_{k}-x_{k})\int_{0}^{1}\frac{df}{d\varsigma_{k}}(\varsigma)dt\right|\\
&\leq&\Big(\sum_{k=1}^{n}|y_{k}-x_{k}|^{2}\Big)^{\frac{1}{2}}\Big[\sum_{k=1}^{n}\Big(\int_{0}^{1}\Big|\frac{\partial
f}{\partial \varsigma_{k}}(\varsigma)\Big|dt\Big)^{2}\Big]^{\frac{1}{2}}\\
&\leq& \sqrt{n}|y-x|\int_{0}^{1}|\nabla f(\varsigma)|dt\\
&\leq&C\sqrt{n}|y-x|\int_{0}^{1}\omega\left(\frac{1}{d(x)-t|y-x|}\right)dt\\
&=&C\sqrt{n}\int_{0}^{|y-x|}\omega\left(\frac{1}{d(x)-t}\right)dt,
\end{eqnarray*}
which implies

\vspace{6pt}

$\ds \frac{1}{|\mathbb{B}^{n}(x,r)|}\int_{\mathbb{B}^{n}(x,r)}|f(\zeta)-f(x)|dV(\zeta)$
\begin{eqnarray*}
&\leq&\frac{C\sqrt{n}}{|\mathbb{B}^{n}(0,r)|}\int_{\mathbb{B}^{n}(0,r)}\left\{\int_{0}^{|\xi|}
\omega\left(\frac{1}{d(x)-t}\right)dt\right\}dV(\xi)\\
&=&\frac{Cn\sqrt{n}}{r^{n}}\int_{0}^{r}\rho^{n-1}\left\{\int_{0}^{\rho}\omega\Big(\frac{1}{d(x)-t}\Big)dt\right\}d\rho\\
&\leq&\frac{Cn\sqrt{n}}{r^{n}}\int_{0}^{r}\left\{\int_{t}^{r}\rho^{n-1}
d\rho\right\}\omega\left(\frac{1}{r-t}\right)dt\\
&\leq&\frac{C\sqrt{n}}{r^{n}}\int_{0}^{r}(r-t)\left(r^{n-1}+r^{n-2}t+\cdots+t^{n-1}\right)\omega\left(\frac{1}{r-t}\right)dt\\
&\leq&\frac{C\sqrt{n}}{r^{n}}r\omega\Big(\frac{1}{r}\Big)\int_{0}^{r}\left(r^{n-1}+r^{n-2}t+\cdots+t^{n-1}\right)dt\\
&=&C\sqrt{n}\left(\sum_{j=1}^{n}\frac{1}{j}\right)r\omega\Big(\frac{1}{r}\Big),
\end{eqnarray*}
where $\varsigma=(\varsigma_{1},\ldots,\varsigma_{n})=yt+(1-t)x$.

Now we   prove that (c)$\Longrightarrow$(a). By Lemma
\ref{lem-g1}, we have
$$|\nabla f(x)|\leq
\frac{n\sqrt{n}}{\rho}\int_{\partial\mathbb{B}^{n}}|f(x+\rho\zeta)-f(x)|d\sigma(\zeta),$$
where $\rho\in(0,d(x)]$. Let $r=d(x)$. Then we have
$$\int_{0}^{r}|\nabla f(x)|\rho^{n} d\rho\leq\sqrt{n}\int_{0}^{r}\Big(n\rho^{n-1}\int_{\partial\mathbb{B}^{n}}|f(x)-f(x+\rho
\zeta)|d\sigma(\zeta)\Big)d\rho,$$ which implies
\begin{eqnarray*}
|\nabla
f(x)|&\leq&\frac{(n+1)\sqrt{n}}{2r^{n+1}}\int_{0}^{r}\Big(n\rho^{n-1}\int_{\partial\mathbb{B}^{n}}|f(x)-f(x+\rho
\zeta)|d\sigma(\zeta)\Big)d\rho\\
&=&\frac{(n+1)\sqrt{n}}{2r|\mathbb{B}^{n}(x,r)|}\int_{\mathbb{B}^{n}(x,r)}|f(\xi)-f(x)|dV(\xi)\\
&\leq&\frac{(n+1)\sqrt{n}C}{2}\omega\Big(\frac{1}{r}\Big)\\
&=&\frac{(n+1)\sqrt{n}C}{2}\omega\left(\frac{1}{d(x)}\right).
\end{eqnarray*}
Therefore, (a)$\Longleftrightarrow$(c). Since
(a)$\Longleftrightarrow$(b) and (a)$\Longleftrightarrow$(c), we
conclude that $$\mbox{(a)}\Longleftrightarrow
\mbox{(b)}\Longleftrightarrow\mbox{(c)}.$$
 The
proof of the theorem is complete. \qed

\section{Landau-Bloch theorem for functions in $\mathcal{H}^{p}(\mathbb{B}^{n}, \mathbb{R}^{n})$ and $\mathcal{PE}_{f}^M$}\label{csw-sec3}

The following lemmas are crucial for the proof of Theorem
\ref{thm1}.

 The following result is a Schwarz-Pick type lemma for harmonic
functions in $\mathcal{H}(\mathbb{B}^{n}, \mathbb{R}^{n})$.

\begin{lem}\label{thm7}
Let  $f\in\mathcal{H}(\mathbb{B}^{n},
\mathbb{R}^{n})$  with $|f(x)|\leq M$ in $\mathbb{B}^{n}$, where $M$
is a positive constant. Then
\be\label{eqb3}
\left|f(x)-\frac{1-|x|}{(1+|x|)^{n-1}}f(0)\right|\leq
M\left[1-\frac{1-|x|}{(1+|x|)^{n-1}}\right].
\ee
\end{lem}
\bpf Without loss of generality, we assume that $f$ is also harmonic
on $\partial\mathbb{B}^{n}.$ We first prove the inequality (\ref{eqb3}). By the
Poisson integral formula, we have
\be\label{xx-1}
f(x)=\int_{\partial\mathbb{B}^{n}}\frac{1-|x|^{2}}{|x-\zeta|^{n}}f(\zeta)\,d\sigma(\zeta),
\ee
where $d\sigma$ denotes the normalized surface measure on $\partial\mathbb{B}^{n}$. By calculations, we have
\begin{eqnarray*}
\left|f(x)-\frac{1-|x|}{(1+|x|)^{n-1}}f(0)\right|&=&\left|\int_{\partial\mathbb{B}^{n}}
\left [\frac{1-|x|^{2}}{|x-\zeta|^{n}}-\frac{1-|x|}{(1+|x|)^{n-1}}\right ]f(\zeta)\,d\sigma(\zeta)\right|\\
&\leq&\int_{\partial\mathbb{B}^{n}}
\left[\frac{1-|x|^{2}}{|x-\zeta|^{n}}-\frac{1-|x|}{(1+|x|)^{n-1}}\right]|f(\zeta)|\,d\sigma(\zeta)\\
&\leq&M\int_{\partial\mathbb{B}^{n}}
\left[\frac{1-|x|^{2}}{|x-\zeta|^{n}}-\frac{1-|x|}{(1+|x|)^{n-1}}\right]\,d\sigma(\zeta)\\
&\leq&M\left[1-\frac{1-|x|}{(1+|x|)^{n-1}}\right]
\end{eqnarray*}
and the proof is complete. \epf


A matrix-valued function $A(x)=\big(a_{i,j}(x)\big)_{n\times n}$is
called {\it matrix-valued and real harmonic function} if each of its
entries $a_{i,j}(x)$ is a real harmonic function from an open subset
$\Omega\subset\mathbb{R}^{n}$ into $\mathbb{R}$.

\begin{lem}\label{lem3.1}
Let $A(x)=\big(a_{i,j}(x)\big)_{n\times n}$ be a matrix-valued
harmonic mapping defined on the ball $\mathbb{B}^{n}(0,r)$. If
$A(0)=0$ and $|A(x)|\leq M$ in $\mathbb{B}^{n}(0,r),$ then
$$|A(x)|\leq M\left[1-\frac{r^{2n-2}(r-|x|)}{(r+|x|)^{2n-1}}\right].
$$
\end{lem}\bpf
For an arbitrary
$\theta=(\theta_1,\ldots,\theta_n)^{T}\in\partial\mathbb{B}^{n}$, we
let
$$ P_{\theta}(x)=A(x)\theta=(p_{1}(x),\ldots,p_{n}(x)).
$$
For every  $\zeta\in\mathbb{B}^{n}$, let
$F_{\theta}(\zeta)=P_{\theta}(r\zeta).$ By Lemma \ref{thm7}, we see
that for all $\zeta\in\mathbb{B}^{n}$,
$$\left|F_{\theta}(\zeta)-\frac{1-|\zeta|}{(1+|\zeta|)^{n-1}}F_{\theta}(0)\right|\leq
M\left[1-\frac{1-|\zeta|}{(1+|\zeta|)^{n-1}}\right],
$$
which gives
$$|P_{\theta}(x)|\leq M\left[1-\frac{r^{n-2}(r-|x|)}{(r+|x|)^{n-1}}\right], ~\mbox{ $|x|<r$}.
$$
The arbitrariness of $\theta$ yields the desired inequality. \epf


\begin{lem}\label{lem3.2}
Let $f\in\mathcal{H}(\mathbb{B}^{n}, \mathbb{R}^{n})$ with $|f(x)|\leq M$ for $x\in\mathbb{B}^{n}$, where
$M$ is a positive constant. Then
$$|f'(x)|\leq M\frac{2|x|+n(1+|x|)}{1-|x|^{2}}.
$$
\end{lem}
\bpf Let $f=(f_{1},\ldots,f_{n})$ and
$\theta=(\theta_{1},\ldots,\theta_{n}) \in\partial\mathbb{B}^{n}$.
Without loss of generality, we assume that $f$ is also harmonic on
$\partial\mathbb{B}^{n}.$ By the Poisson integral formula, we find that
$$f(x)=\int_{\partial\mathbb{B}^{n}}\frac{1-|x|^{2}}{|x-\zeta|^{n}}f(\zeta)\,d\sigma(\zeta),
$$
where $d\sigma$ denotes the normalized surface measure on $\partial\mathbb{B}^{n}$. Clearly,
\be\label{eq-1}
\int_{\partial\mathbb{B}^{n}}\frac{d\sigma(\zeta)}{|x-\zeta|^{n}} =\frac{1}{1-|x|^{2}}.
\ee
For each$j,k\in\{1,\ldots,n\}$, we have
$$(f_{j}(x))_{x_{k}}=\int_{\partial\mathbb{B}^{n}}
\frac{-2x_{k}|x-\zeta|^{2}-n(1-|x|^{2})(x_{k}-\zeta_{k})}{|x-\zeta|^{n+2}}f_{j}(\zeta)
\, d\sigma(\zeta),
$$
which gives
\begin{eqnarray*}
\left |\sum_{k=1}^{n}(f_{j}(x))_{x_{k}}\cdot\theta_{k}\right |^{2}&=&
\left |\sum_{k=1}^{n}\int_{\partial\mathbb{B}^{n}}
\frac{[2x_{k}|x-\zeta|^{2}+n(1-|x|^{2})(x_{k}-\zeta_{k})]\theta_{k}}
{|x-\zeta|^{n+2}}f_{j}(\zeta)\, d\sigma(\zeta)\right |^{2}\\
&=&\left |\int_{\partial\mathbb{B}^{n}}
\frac{\sum_{k=1}^{n}[2x_{k}|x-\zeta|^{2}+n(1-|x|^{2})(x_{k}-\zeta_{k})]\theta_{k}}
{|x-\zeta|^{n+2}}f_{j}(\zeta)\, d\sigma(\zeta)\right |^{2}\\
&\leq&\left [\int_{\partial\mathbb{B}^{n}}\frac{[2|x|\,|x-\zeta|^{2}+n(1-|x|^{2})|x-\zeta|]|f_{j}(\zeta)|}
{|x-\zeta|^{n+2}}\,d\sigma(\zeta)\right ]^{2}\\
&\leq&\left [\int_{\partial\mathbb{B}^{n}}\frac{[2|x|\,|x-\zeta|+n(1-|x|^{2})]^{2}}
{|x-\zeta|^{n+2}}\,d\sigma(\zeta)\right ]\\
&&\times \left [\int_{\partial\mathbb{B}^{n}}\frac{|f_{j}(\zeta)|^{2}}{|x-\zeta|^{n}}\,d\sigma(\zeta)\right ]
\end{eqnarray*}
Then the relation \eqref{eq-1} shows
\begin{align*}
\sum_{j=1}^{n}\left |\sum_{k=1}^{n}(f_{j}(x))_{x_{k}}\cdot\theta_{k}\right |^{2}&\leq
\left [\int_{\partial\mathbb{B}^{n}}\frac{[2|x|\,|x-\zeta|+n(1-|x|^{2})]^{2}}
{|x-\zeta|^{n+2}}\,d\sigma(\zeta)\right ]\\
~~~~~~~~~~&\times \left [\int_{\partial\mathbb{B}^{n}}\frac{\sum_{j=1}^{n}|f_{j}(\zeta)|^{2}}
{|x-\zeta|^{n}}\,d\sigma(\zeta)\right ]\\
&\leq\frac{M^{2}}{1-|x|^{2}}\left [\int_{\partial\mathbb{B}^{n}}\frac{[2|x|\,|x-\zeta|+n(1-|x|^{2})]^{2}}
{|x-\zeta|^{n+2}}\,d\sigma(\zeta)\right]\\
&\leq\frac{M^{2}}{1-|x|^{2}}\left [\int_{\partial\mathbb{B}^{n}}
\frac{[2|x|+n(1+|x|)]^{2}}{|x-\zeta|^{n}}\,d\sigma(\zeta)\right]\\
&\leq\frac{M^{2}[2|x|+n(1+|x|)]^{2}}{(1-|x|^{2})^{2}}
\end{align*}
whence
$$|f'(x)|\leq M\frac{2|x|+n(1+|x|)}{1-|x|^{2}}.
$$
The proof of this lemma is complete. \epf

\begin{Lem}{\rm (\cite[Lemma 4]{LX})}\label{LemA}
Let $A$ be an $n \times n$ real $($or complex$)$ matrix with $|A|\neq0$. Then for any unit
vector $\theta\in\partial \mathbb{B}^{n}$, the inequality
$$|A\theta|\geq\frac{|\det A|}{|A|^{n-1}}
$$
holds.
\end{Lem}

\subsection*{Proof of Theorem \ref{thm1}} Without loss of generality, we assume that $f$ is also harmonic on
$\partial\mathbb{B}^{n}$, where $n\geq3$. By the Poisson integral
representation, we have
$$f(x)=\int_{\partial\mathbb{B}^{n}}\frac{1-|x|^{2}}{|x-\zeta|^{n}}f(\zeta)\,d\sigma(\zeta)
$$
in $\mathbb{B}^{n}$. By Jensen's inequalities, we obtain
$$|f(x)|^{p}\leq\int_{\partial\mathbb{B}^{n}}\frac{1-|x|^{2}}{|x-\zeta|^{n}}
|f(\zeta)|^{p}\,d\sigma(\zeta)\leq\frac{2\|f\|_{p}^{p}}{(1-|x|)^{n-1}}
$$
which gives
$$|f(x)|\leq\frac{2^{1/p}K_{0}}{(1-|x|)^{(n-1)/p}},
$$
where $K_{0}=\|f\|_{p}$. For $\zeta\in\mathbb{B}^{n}$ and for a fixed $r\in(0,1)$, let
$F(\zeta)=f(r\zeta)/r.$ Then
$$|F(\zeta)|\leq\frac{2^{1/p}K_{0}}{r(1-r)^{(n-1)/p}}=K(r).
$$
For each $\zeta\in\mathbb{B}^{n}(0,\sqrt{2}/2)$, using Lemma
\ref{lem3.2}, we have
\begin{eqnarray*}
|F'(\zeta)-F'(0)|
&\leq&|F'(0)|+|F'(\zeta)|\\
&\leq&nK(r)+\frac{K(r)[n+(n+2)|\zeta|]}{1-|\zeta|^{2}}\\
&\leq&K(r)[(3+\sqrt{2})n+2\sqrt{2}],
\end{eqnarray*} which implies $F'(\zeta)-F'(0)$ is a bounded  matrix-valued and real harmonic
function in $\mathbb{B}^{n}(0,\sqrt{2}/2).$ By Lemma \ref{lem3.1},
for each $\zeta\in\mathbb{B}^{n}(0,\sqrt{2}/2)$, we have
\begin{eqnarray*}
|F'(\zeta)-F'(0)|&\leq&
M(r)\left[1-\frac{\left(\frac{\sqrt{2}}{2}\right)^{n-2}(\frac{\sqrt{2}}{2}-|\zeta|)}{(\frac{\sqrt{2}}{2}+|\zeta|)^{n-1}}\right]\\
&\leq&M(r)\cdot\frac{C_{n-1}^{1}\left(\frac{\sqrt{2}}{2}\right)^{n-2}|\zeta|+\cdots+C_{n-1}^{n-1}|\zeta|^{n-1}+
\left(\frac{\sqrt{2}}{2}\right)^{n-2}|\zeta|}{(\frac{\sqrt{2}}{2}+|\zeta|)^{n-1}}\\
&\leq&M(r)|\zeta|\frac{\left[(1+\frac{\sqrt{2}}{2})^{n-1}+\Big(\frac{\sqrt{2}}{2}\Big)^{n-2}-\Big(\frac{\sqrt{2}}{2}\Big)^{n-1}\right]}
{(\frac{\sqrt{2}}{2}+|\zeta|)^{n-1}}\\
&\leq&M(r)\left[(1+\sqrt{2})^{n-1}+\sqrt{2}-1\right]|\zeta|,
\end{eqnarray*}
where $M(r)=K(r)[(3+\sqrt{2})n+2\sqrt{2}]$   and $C_{n}^{k}={n\choose
k}$ ($k=1,2, \ldots, n$) denote the binomial coefficients.

 Since for each $\theta\in\partial\mathbb{B}^{n}$,
Lemmas \ref{lem3.2} and \Ref{LemA} imply
$$|F'(0)\theta|\geq\frac{J_{F}(0)}{|F'(0)|^{n-1}}\geq \frac{1}{\big [nK(r)\big ]^{n-1}}.
$$

Let $\zeta'$ and $\zeta''$ be two distinct points in
$\mathbb{B}^{n}(0,\rho(r))$ with
$$\rho(r)=\frac{1}{[nK(r)]^{n-1}M(r)[(1+\sqrt{2})^{n-1}+\sqrt{2}-1]},
$$
and let $[\zeta',\zeta'']$ denote the segment connecting $\zeta'$ and $\zeta''$. Set
$$ d\zeta=\left(\begin{array}{cccc}
d\zeta_{1}   \\
\vdots \\
d\zeta_{n}
\end{array}\right).
$$
Then we have
\begin{eqnarray*}
|F(\zeta')-F(\zeta'')|&\geq&
\left|\int_{[\zeta',\zeta'']}F'(0)d\zeta\right|
-\left|\int_{[\zeta',\zeta'']}(F'(\zeta)-F'(0))\,d\zeta\right|\\
&>&|\zeta'-\zeta''|\left\{\frac{1}{[nK(r)]^{n-1}}-M(r)\left[(1+\sqrt{2})^{n-1}+\sqrt{2}-1\right]\rho(r)\right\}\\
&=&0.
\end{eqnarray*}
This observation shows that $F$ is univalent in $\mathbb{B}^{n}(0,\rho(r))$. Furthermore, for any
$\zeta_{0}$ with $|\zeta_{0}|=\rho(r)$, we have
\begin{eqnarray*}
|F(\zeta_{0})-F(0)|&\geq&
\left|\int_{[0,\zeta_{0}]}F'(0)d\zeta\right|
-\left|\int_{[0,\zeta_{0}]}(F'(\zeta)-F'(0))\, d\zeta\right|\\
&\geq&\rho(r)\left\{\frac{1}{[nK(r)]^{n-1}}-M(r)\left[(1+\sqrt{2})^{n-1}+\sqrt{2}-1\right]\rho(r)/2\right\}\\
&=&\frac{\rho(r)}{2[nK(r)]^{n-1}}\\
&>&0.
\end{eqnarray*}
Therefore, $f(\mathbb{B}^{n})$ contains a univalent ball
$\mathbb{B}^{n}(0,R)$, where
\begin{eqnarray*}
R&\geq&\max_{0<r<1}\left\{\frac{\rho(r)}{2[nK(r)]^{n-1}}\right\}\\
&=&\max_{0<r<1}\left\{\frac{1}{2[nK(r)]^{2n-2}M(r)[(1+\sqrt{2})^{n-1}+\sqrt{2}-1]}\right\}.
\end{eqnarray*}
 The
theorem is proved. \qed


\subsection*{Proof of Theorem \ref{thm6-x}} If we suppose that this result is not true, then
there is a sequence $\{a_k\}$    and a sequence of functions
$\{u_k\}$ with $u_{k}\in\mathcal{PE}_{f}^{M}$, such that $\{a_k\}$
tends to $0$ and $a_k \notin u_k(\mathbb{B}^{n})$, where $a_k>0$ for
$k\in\{1,2,\ldots\}$. By \cite[Theorem 4.6 and Corollary 4.7]{gt},
we know that there is a subsequence $\{g_{k}\}$ of $\{u_{k}\}$ which
converges uniformly on compact subsets of $\mathbb{B}^{n}$ to a
function $g $. Note that for each $k$, the function $h_k=g_{k}-g_1$ is harmonic.
Hence the sequence $\{h_k\}$ converges uniformly on compact subsets of
$\mathbb{B}^{n}$ to  $ g-g_1$  and therefore, the partial derivatives
of $g_k$ converge uniformly on compact subsets of $\mathbb{B}^{n}$ to
the partial derivatives of $ g$. In particular,  ${g_k} (0)
\rightarrow g(0)$ and $J_{g_k} (0) \rightarrow J_{g} (0)$, and
therefore, $g \in \mathcal{PE}_{f}^M$. Since $J_{g} (0)-1=|g(0)|=0,$
there are $0 < r_0 <1$ and $c_1
> 0$  such that  $J_{g} >0$  on  $
\overline{\mathbb{B}^{n}(0,r_{0})}$, $g(\mathbb{B}^{n}(0,r_{0}))
\supset \overline{\mathbb{B}^{n}(0,c_1)}$  and $|g(x)| \geq c_1$ for
$x \in
\partial\mathbb{B}^{n}(0,r_{0})$.

Set $c_2=c_1/2$, $B_{r_{0}}=\mathbb{B}^{n}(0,r_{0})$ and $B_{c_{2}}=
\mathbb{B}^{n}(0,c_2) $. Then there is a $k_0$ such that $|g_{k}(x)|
\geq c_2$ for $k \geq k_0$  and  $J_{g_k} >0$  on
$\overline{B_{r_{0}}}$. Since $\deg(g_k,B_{r_{0}},0)\geq 1$, by the
degree property (II) in page \pageref{(II)},   we see  that
$\deg(g_k,B_{r_{0}},y)\geq 1$ for $y \in B_{c_{2}}$ and $k \geq
k_0$. Hence  $g_k(B_{r_{0}})\supset B_{c_{2}}$ for $k \geq k_0$ and
this leads a contradiction. The proof of the theorem is complete.
\qed

\end{document}